\documentclass[12pt]{article}
\usepackage{amsfonts,amssymb,amsmath}
\usepackage{graphicx,graphics,psfrag}
\usepackage{amsmath}
\DeclareMathOperator\arctanh{arctanh}

\usepackage[english]{babel}
\makeatletter
\@addtoreset{equation}{section}
\makeatother

\newcounter{enunciato}[section]

\newtheorem{ittheorem}{Theorem}
\newtheorem{itlemma}{Lemma}
\newtheorem{itproposition}{Proposition}
\newtheorem{itdefinition}{Definition}
\newtheorem{itcorollary}{Corollary}
\newtheorem{itconjecture}{Conjecture}

\newenvironment{theorem}{\addtocounter{enunciato}{1}
\begin{ittheorem}}{\end{ittheorem}}

\newenvironment{corollary}{\addtocounter{enunciato}{1}
\begin{itcorollary}}{\end{itcorollary}}

\newcommand{\halmos}{\rule{1ex}{1.4ex}}

\newenvironment{proof}{\noindent {\em Proof}.\,\,}
{\hspace*{\fill}$\halmos$\medskip}

\def \ba {\begin{array}}
\def \ea {\end{array}}

\def \N {{\mathbb N}}

\begin{document}
\title{Phase Transition in Conditional Curie-Weiss Model}

\author{
\renewcommand{\thefootnote}{\arabic{footnote}}
Alex Akwasi \ Opoku \footnotemark[1]
, Kwame\ Owusu Edusei \footnotemark[2]
\renewcommand{\thefootnote}{\arabic{footnote}}
  and Richard Kwame Ansah
\footnotemark[3]
}

\footnotetext[1]{
alex.opoku@uenr.edu.gh
}
\footnotetext[2]{
92sherlock@gmail.com
}
\footnotetext[3]{
richard.ansah@uenr.edu.gh
}
\footnotetext[1]{
$^{,2,3}$Mathematics and Statistics Department, 
University of Energy and Natural Resources,
P. O. Box  214, Sunyani, Ghana.
}

%

\maketitle

\begin{abstract}
This paper  proposes a conditional Curie-Weiss    model as a model for opinion formation  in a  
 society polarized along two opinions, say opinions 1 and 2. The model  comes with interaction strength $\beta>0$ and 
 bais $h$. Here the population in question is divided into three main groups, namely:
\begin{enumerate}
\item Group one consisting of individuals who have decided on opinion 1. Let the  proportion of this
 group be given by $s$.
\item Group two consisting of individauls who have chosen opinion 2. Let $r$ be their 
proportion. 

\item Group three consisting of individuals  who are yet to decide and they will decide based on 
their environmental conditions.  Let $1-s-r$ be the proportion of this group.
\end{enumerate}
We show that the specific magnetization of the  associated conditional Curie-Weiss model  
has a  first order phase transition (discontinuous jump in specific magnetization) at 
$\beta^*=\left(1-s-r\right)^{-1}$. It is also shown that not all the discontinuous jumps in 
magnetization will result in    phase change. We point out how an extention of this model could 
serve as a random field Curie-Weiss model where the random field distribution has  nonvanishing mean.

\vskip .2in
\medskip\noindent
{\it AMS} 2000 {\it subject classifications.} 60F10, 60K37, 82B27.

\vskip .2in
\medskip\noindent
{\it Acknowledgment.}
The authors thank  the staff at the Mathematics Department of University of Energy and Natural
 Resources for their support and kindness during the period this paper was written.

\end{abstract}


\section{Introduction and main results}
\label{S1}

Phase transition in a physical system corresponds to a change from one phase (behaviour) of 
the  physical system to another through changes in the  environment of the system. For 
instance the liquid-vapour transition of water when it is boiled. In the words of  J. Willard 
Gibbs a phase transition is a singularity in thermodynamic behavior, i.e. singularity in the 
free energy  of the system \cite{Kad2009}. The nature of the singularity, according to  
P. Ehrenfest, determines the order of the phase transition. For instance,  first order phase 
transition  is associated with a discontinous jump in some thermodynamic quantity called the 
order parameter of the system. An important theory that has offered a lot of insight into the 
study  of phase transitions is the mean field theory. This theory says that each component of 
the physical systems  feels the average influence of all the other components of the system. 
Johannes Diderik van der Waals was the first to derive  mean field theory in the 1870s in his  
attempt at understanding the 1869  experimental data of  T. Andrews \cite{Andrew89} that 
showed a phase transition curve separating the  liquid-vapour phases of  fluids. In 1895,  
Pierre Curie observed that a ferromagnet  admits a behaviour similar to that of fluids \cite{PC95}. 
A mean field theory for ferromagnet was derived in 1907 by Pierre Weiss \cite{PW07}. In the late 
1960s Mark Kac developed a model for ferromagnets where every magnetic moment interacts with 
every other magnetic moment and this model is known in the literature as the  Curie-Weiss (CW) 
model \cite{Kac68}.

The Curie-Weiss model after its intoduction has found other applications apart from what it was originally 
designed to do. Notable among these applications are formation of opinion in societies \cite{BP98, BP98B},  
immigrants' integration \cite{ABC2014, BCS2014, CV05}, democratization \cite{CPRS06}, etc. 

The present paper studies phase transition in Curie-Weiss model conditional on having certain minimum 
proportions of magnetic moments following each of the two possible spin alignments. This could serve as a model 
for opinion formation in a society segregated along two opinions. Here the minimum proportions are proportions of
individuals who have fixed their orientation on one of the two opinions and they  will never change their views. 
Phase transition here will imply emergence of consensus and lack of it  will imply wild flactuation 
of  group opinion  and  never settling on a specific collective opinion or decision, i.e. non-consensus 
emergence \cite{BP98}.

In a followup paper, we will study the case where the minimum proportions are fixed according to some 
distribution and we ask for the effect of the disorder in these proportions on the phase transtions in the 
associated quenched and annealed models. This will 
provide a natural example of random field Ising type of model \cite{AB2006} with non centered 
random field. The present paper is the first step towards investigating these class of random 
field Ising models.  Further, such models will naturally set the stage for  studying spin models 
on  site percolation clusters generated from  a 
random process of assigning three different colours to the vertices of the underlying  graph. 
The spins on one of the three clusters, generated from the three 
colours, will be fixed to +1, one of the remaining two will also be set to -1 and the spins 
on the  remaining cluster could pick any of the spin values. We then ask the question of 
phase transition in Ising spin model on such a decorated graph.

The rest of the paper is organized as follows: In Section \ref{S1.1} we recall the Curie-Weiss model and 
collect some fact about it. Section \ref{S1.2} is devoted to  defining our conditional Curie-Weiss model. 
The main results of the paper are collected in Section \ref{S1.3}. The results in Section \ref{S1.3} are 
discussed in Section \ref{S2} and the proofs of the main results are found in Section \ref{S3}.


\subsection{The Curie-Weiss Model}
\label{S1.1}

Let $N$ be a positive integer, $V_N=\left\{1, 2, \ldots, N \right\}$ and $E_N$ be the vertex set  
and the edge set respectively for a complete graph with $N$ vertices. Denote by $\Omega_N=\{ -1,+1\}^N$, the set of 
configurations of the system indexed by the elements of  $V_N$.
The Curie-Weiss model is a probability measure $\mu_{N}$ on $\Omega_N$ given by 
\begin{equation}\label{CWmodel}
\mu_N(\sigma)=\frac{1}{Z_N}\exp\left( -\beta H_N(\sigma)\right), 
\end{equation}
where for any $\sigma \in \Omega_N$
\begin{equation}\label{Hamil}
\begin{split}
H_N(\sigma)&=-\frac{1}{2N}\sum_{(i,j)\in E_N}\sigma_i\sigma_j-h\sum_{i=1}^N\sigma_i \cr
&=-N\left[\frac{1}{2}\left(\frac{1}{N}\sum_{i=1}^N\sigma_i \right)^2-\frac h N\sum_{i=1}^N
\sigma_i\right]-\frac12\cr
&= -N\left[\frac{1}{2} m_N^2+hm_N\right]-\frac12.
\end{split}
\end{equation}
Here $h$ is a real number and  $\beta$ is a positive real number. The quantity $Z_N$ is a normalization 
term called the partition function of the model. $h$ and $\beta$ are the parameters of the model
called the interaction bais and interaction strength respectively. $h$ is usually called an external 
field and $\beta$ is the inverse temperature.
The function $H_N$ on configurations set $\Omega_N$ is called the Hamiltonian/ energy function 
of the model. The first term in the Hamiltonian turns to align pairs of spins, while the second turns to align 
spins in the direction of the external field $h$. It is known in the literature that the Curie-Weiss models 
undergoes a phase transition, namely the  specific magnetization
\begin{equation}
m(\beta,h)=\lim_{N\rightarrow\infty}\frac{1}{N}\int_{\Omega_N} \left(\sum_{i=1}^N\sigma_i\right)
d\mu_N(\sigma)
\end{equation}
exhibits a discontinuity as $h\rightarrow 0$ for $\beta>1$. In fact, it  is proved in  Theorem IV.4.1 (b) 
of \cite{Ellis1985} that:
\begin{theorem}\label{mainfact}
For $\beta>0$, $h\neq 0$ and $0<\beta\leq 1$, $h=0$, the specific magnetization $m(\beta,h)$ is equal
to $z(\beta,h)$, where $z(\beta,h)$ is the minimizer of the function 
\begin{equation}
i_{\beta,h}(z)=-\frac12\beta z^2-\beta h z+\frac{1-z}{2}\log(1-z)+\frac{1+z}{2}\log(1+z).
\end{equation}
In particular, 
\begin{equation}
m(\beta,\pm)=\lim_{h\rightarrow 0^{\pm}}m(\beta,h)=\left\{\begin{array}{ll}
z(\beta, 0)=0 & \mbox{for $0<\beta\leq 1$},\\
z(\beta, \pm)\neq 0 &\mbox{ for $\beta>1$}.
\end{array}\right.
\end{equation}

\end{theorem}
Thus when $\beta>1$,  any small change in $h$, no matter how small it is, that we make around  $h=0$
  will lead to a dramatic change in the behaviour of the system. 

\subsection{Conditional Curie-Weiss Model}
\label{S1.2}
 For any positive integer $N$, let $s_N$, $r_N$ and $t_N$ be positive real numbers such that 
$r_N+s_N+t_N=1$, $r_N\rightarrow r$,  $s_N\rightarrow s$ and $t_N\rightarrow t$ as $N\rightarrow\infty.$ 
Further, let $s_N$, $r_N$ and $t_N$ be such that there is a partition 
$V_{N,r_N}$, $V_{N,s_N}$ and $V_{N,t_N}$  of $V_N$  such that $r_N=\frac{|V_{N,r_N}|}{N}$, 
$s_N=\frac{|V_{N,s_N}|}{N}$ and $t_N=\frac{|V_{N,t_N}|}{N}$. 
Here we write  $|A|$ for the cardinality of a set $A$. Define a subset $\Omega_{N,s,r}$ of 
$\Omega_N$ consisting of configurations $\sigma=(\sigma_i)_{1\leq i\leq N}$ such  that 
\begin{equation}\label{configsr}
\sigma_i=\left\{ \begin{array}{ll}
1 & \mbox{ if  $i\in V_{N, s_N}$};\\
-1 &\mbox{ if  $i\in V_{N,r_N}$};\\
\eta_{i} & \mbox{ if  $i\in V_{N,t_N}$},

\end{array}\right.
\end{equation}
where $\eta_i\in\{+1, -1\}$. Thus all sites labelled by $V_{N,s_N}$ have fixed spin value  of $+1$, those on 
$V_{N,r_N}$ are  fixed to the spin value  of $-1$ and those labelled by $V_{N,t_N}$ are free 
to choose any of the two spin values.  Here $+1$  corresponds to \emph{opinion one}  and $-1$ to \emph{opinion two}. 
Thus  the part of the population indexed by $V_{N, s_N}$ are those who have decided on opinion one 
and those with label set $V_{N,r_N}$ are those who have chosen opinion two with those in 
$V_{N,t_N}$ yet to decide.
In this paper we study a conditional distribution $\mu_{N,s,r}$ of the Curie-Weiss model  $\mu_N$  
\eqref{CWmodel} conditional  on spins in $\Omega_{N,s,r}$. Thus $\mu_{N,s,r}$ is given by 
\begin{equation}\label{condCWmodel}
\begin{split}
\mu_{N,s,r}(\sigma)&=\mu_{N}(\sigma|\Omega_{N,s,r})\cr
&=\frac{\mu_N(\sigma)}{\mu_N(\Omega_{N,s,r})}\cr
&=\frac{e^{-\beta H_N(\sigma)}}{\sum_{\eta\in\Omega_{N,s,r}} e^{-\beta H_N(\eta)}}\cr
&=\frac{e^{-\beta H_N(\sigma)}}{\widetilde{Z}_{N,s,r}}\
\end{split}
\end{equation}
Thus $\mu_{N,s,r}$ is a probability measure on $\Omega_{N,s,r}$.
The above conditional probability is well defined as for all $N\in\N$
\begin{equation}
\begin{split}
\mu_N(\Omega_{N,s,r})&=\dfrac{\sum_{\sigma\in\Omega_{N,s,r}}e^{-\beta H_N(\sigma)}}{
\sum_{\tilde\sigma\in\Omega_{N}}e^{-\beta H_N(\tilde\sigma)}}\cr
&\geq 2^{(|V_{N,t_N}|-N)} 
e^{-\beta\sup_{\sigma,\tilde\sigma\in\Omega_N}\left|H_N(\sigma)-H_N(\tilde{\sigma})\right|}\cr
&\geq 2^{(|V_{N,t_N}|-N)} \,\,
e^{-\beta\left( \frac{N-1}{2}+2N|h|\right)}\cr
&>0.
\end{split}
\end{equation}
We are ready to state the main results of the paper in the next subsection. Before we do this, let us define 
 \begin{equation}
 m(\beta,s,r,h)=\lim_{N\rightarrow\infty}\frac1N\int_{\Omega_{N,s,r}}\left(
 \sum_{i=1}^N\sigma_i\right)d\mu_{N,s,r}(\sigma).
 \end{equation}
 \newpage
 \subsection{Main Results}\label{S1.3}
 In what follows, we always assume $s,\; r\geq 0$ and  $s+r<1$.

\begin{figure}[htbp]
\vspace{-1.2cm}
\begin{center}
\includegraphics[scale = 0.45]{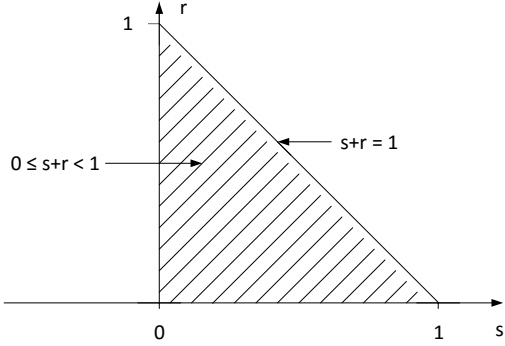}
\end{center}
\vspace{-7.0cm}
\caption{\small \emph{ The shaded region in the $sr$-- plane is the region of interest for this study.}}
\label{fig3}
\vspace{0cm}
\label{fig}
\end{figure}

 \begin{theorem}\label{maintheorem}
For $\beta>0$, $h\neq r-s$ and $0<\beta\leq (1-s-r)^{-1}$, $h=r-s$, the specific magnetization 
$m(\beta,s,r,h)$ is given by
\begin{equation}
m(\beta,s,r,h)=s-r+(1-s-r) z(\beta,s,r,h), 
\end{equation}
 where $z(\beta,s,r,h)$ is the minimizer of the function 
\begin{equation}
i_{\beta,s,r,h}(z)=-\frac12\beta (1-s-r) z^2-\beta(s-r+ h) z+\frac{1-z}{2}\log(1-z)+\frac{1+z}{2}\log(1+z), 
\end{equation}
for $-1\leq z\leq 1.$ In particular, 
\begin{equation}\label{meq}
\begin{split}
&m(\beta,s,r,(r-s)^\pm)=\lim_{h\rightarrow (r-s)^{\pm}}m(\beta,s,r,h)\cr
&=\left\{\begin{array}{ll}
s-r, & \mbox{for $0<\beta\leq (1-s-r)^{-1}$},\\\\
s-r+(1-s-r)z(\beta,s,r,(r-s)^\pm)\neq s-r, &\mbox{ for $\beta>(1-s-r)^{-1}$},
\end{array}\right.
\end{split}
\end{equation}
where 
\begin{equation}
\begin{split}
z(\beta,s,r,(r-s)^\pm)=\lim_{h\rightarrow (r-s)^{\pm}} z(\beta,s,r,h).
\end{split}
\end{equation}

\end{theorem}
The following corollaries list what happens to the limit in \eqref{meq} as we move through 
the shaded region of Figure \ref{fig}.
\begin{corollary}\label{cor1}
Suppose $s=r$ and $0\leq s <\frac12$. Then

\begin{equation}
\begin{split}
m(\beta,s,s,0^\pm)&=\lim_{h\rightarrow 0^{\pm}}m(\beta,s,s,h)\cr
&=\left\{\begin{array}{ll}
0,& \mbox{for $0<\beta\leq (1-2s)^{-1}$},\\
(1-2s)z(\beta,s,s,0^\pm)\neq 0, &\mbox{ for $\beta>(1-2s)^{-1}$}.
\end{array}\right.
\end{split}
\end{equation}
In particular, if $\beta>(1-2s)^{-1}$ then  there is a discontinuity in $m(\beta,s,r,h)$, as a function of 
$h$, at $h=0$ and this discontinuity is 
followed by change in state, i.e. $m(\beta,s,s,0^+)=-m(\beta,s,s,0^-).$

\end{corollary}

\begin{figure}[htbp]
\vspace{-3.0cm}
\begin{center}
\includegraphics[scale = 0.44]{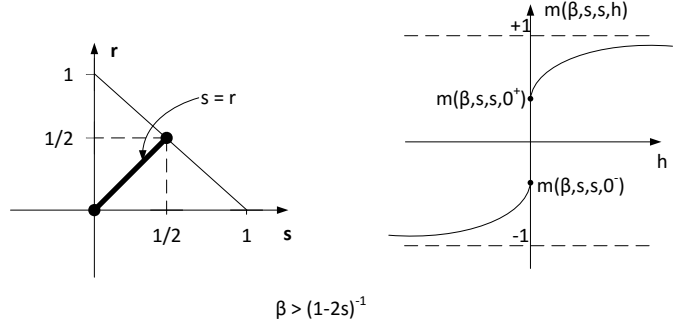}
\end{center}
\vspace{-6.6cm}
\caption{\small \emph{ Discontinuity in the map $h\mapsto m(\beta,s,r,h)$ at $h=0$ for $\beta>
 (1-2s)^{-1}$ and $s=r$.}}
\label{fig3}
\vspace{0cm}
\end{figure}

\begin{corollary}\label{cor2}
Suppose $\frac12\leq s<1$ or $\frac12\leq r<1$ and $\frac12\leq s+r<1$. Then
\begin{equation}
\begin{split}
&m(\beta,s,r,(r-s)^\pm)=\lim_{h\rightarrow (r-s)^{\pm}}m(\beta,s,r,h)\cr
&=\left\{\begin{array}{ll}
s-r, & \mbox{for $0<\beta\leq (1-s-r)^{-1}$},\\
s-r+(1-s-r)z(\beta,s,r,(r-s)^\pm)\neq s-r, &\mbox{ for $\beta>(1-s-r)^{-1}$}.
\end{array}\right.
\end{split}
\end{equation}
Here there is a discontinuity in the map  $h\mapsto m(\beta,s,r,h)$ at $h=r-s$ for $\beta>(1-s-r)^{-1}$. 
This discontitnuity is not followed  by change in state, i.e. $m(\beta,s,r,(r-s)^+)$  and \\$m(\beta,s,r,(r-s)^-)$ 
 are not equal but have the same sign as $s-r$.
\end{corollary}

\begin{figure}[htbp]
\vspace{-5.39cm}
\begin{center}
\centering
\includegraphics[scale = 0.44]{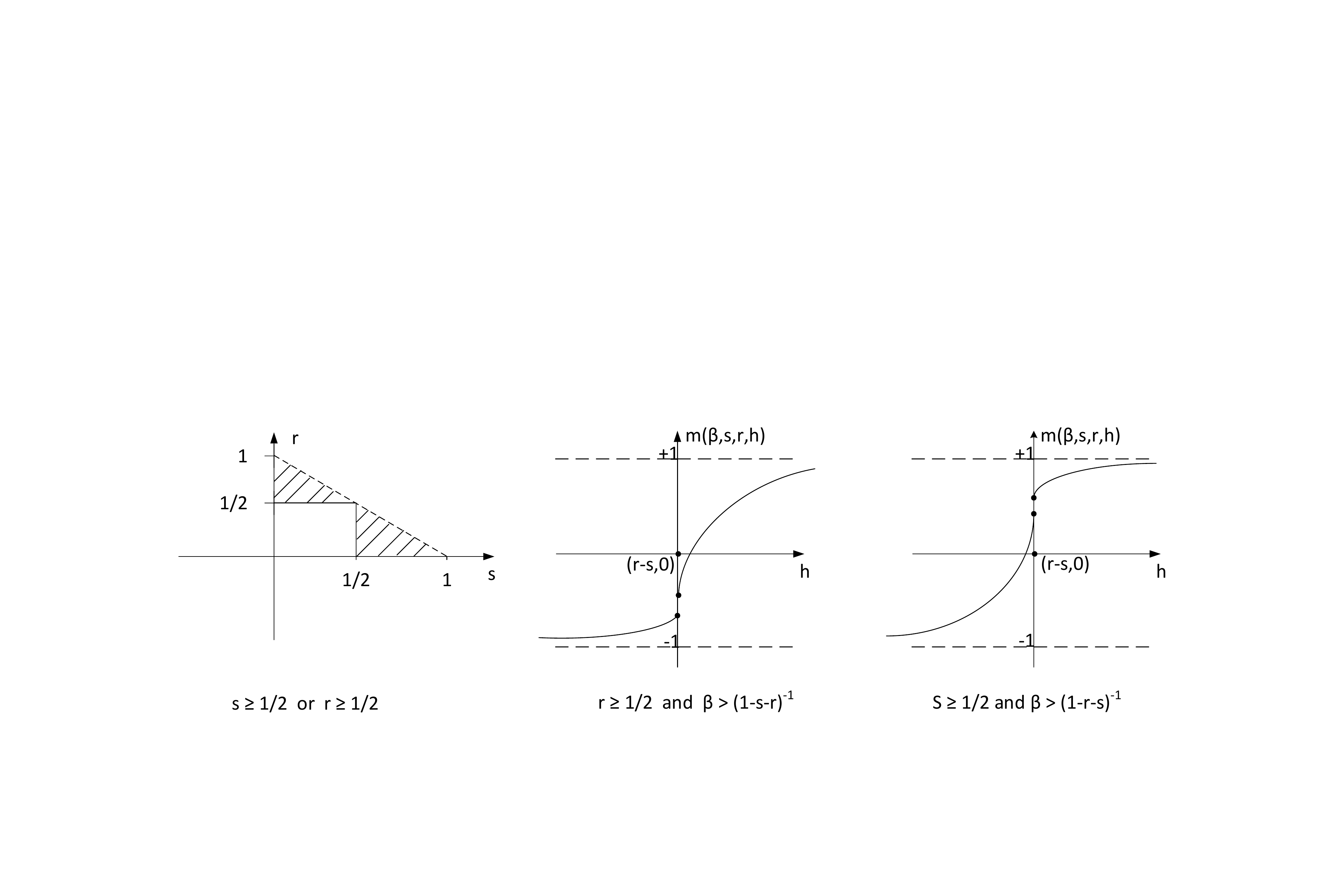}
\end{center}
\vspace{-3.1cm}
\caption{\small \emph{ Discontinuity in the map $h\mapsto m(\beta,s,r,h)$ at $h=r-s$ for $\beta> (1-s-r)^{-1}$
 and for  choices of  $(s,r)$ from the shaded region in the $sr$-plane. The shaded region corresponding 
 to $\frac12\leq r<1$ has $s-r<0$ and that  of  $\frac12\leq s<1$ has $s-r>0$.}}
\label{fig3}
\vspace{0cm}
\end{figure}

\begin{corollary}\label{cor3}
Suppose $0\leq s, r<\frac12$ and $s\neq r$. Then
\begin{equation}
\begin{split}
&m(\beta,s,r,(r-s)^\pm)=\lim_{h\rightarrow (r-s)^{\pm}}m(\beta,s,r,h)\cr
&=\left\{\begin{array}{ll}
s-r, & \mbox{for $0<\beta\leq (1-s-r)^{-1}$},\\
s-r+(1-s-r)z(\beta,s,r,(r-s)^\pm)\neq s-r, &\mbox{ for $\beta>(1-s-r)^{-1}$}.
\end{array}\right.
\end{split}
\end{equation}
There is discontinuity in $h\mapsto m(\beta,s,r,h)$ at $h=r-s$ for $\beta>(1-s-r)^{-1}$. 
Here we have the followed  scenarios:
\begin{enumerate}
\item  The discontitnuity is not followed  by change in state, i.e. $m(\beta,s,r,(r-s)^+)$  and 
\\$m(\beta,s,r,(r-s)^-)$ have the same sign as that of $s-r$ for all 
\begin{equation}
(1-s-r)^{-1}<\beta<\frac{1}{r-s}\arctanh\left(\frac{r-s}{1-s-r}\right).
\end{equation}
\item For 
\begin{equation}
\beta=\frac{1}{r-s}\arctanh\left(\frac{r-s}{1-s-r}\right),
\end{equation}
 The discontinuity is followed 
by change from an ordered state (a plus or minus phase depending on the sign of $s-r$) to a 
disordered state with zero specific magnetization. 
\item Finally, for 
\begin{equation}
\beta>\frac{1}{r-s}\arctanh\left(\frac{r-s}{1-s-r}\right),
\end{equation}
the discontinuity is followed by change in state i.e. $m(\beta,s,r,(r-s)^+)>0$  and 
$m(\beta,s,r,(r-s)^-)<0$.

\end{enumerate}

\end{corollary}

\begin{figure}[htbp]
\vspace{-2.29cm}
\begin{flushleft}
\includegraphics[scale = 0.45]{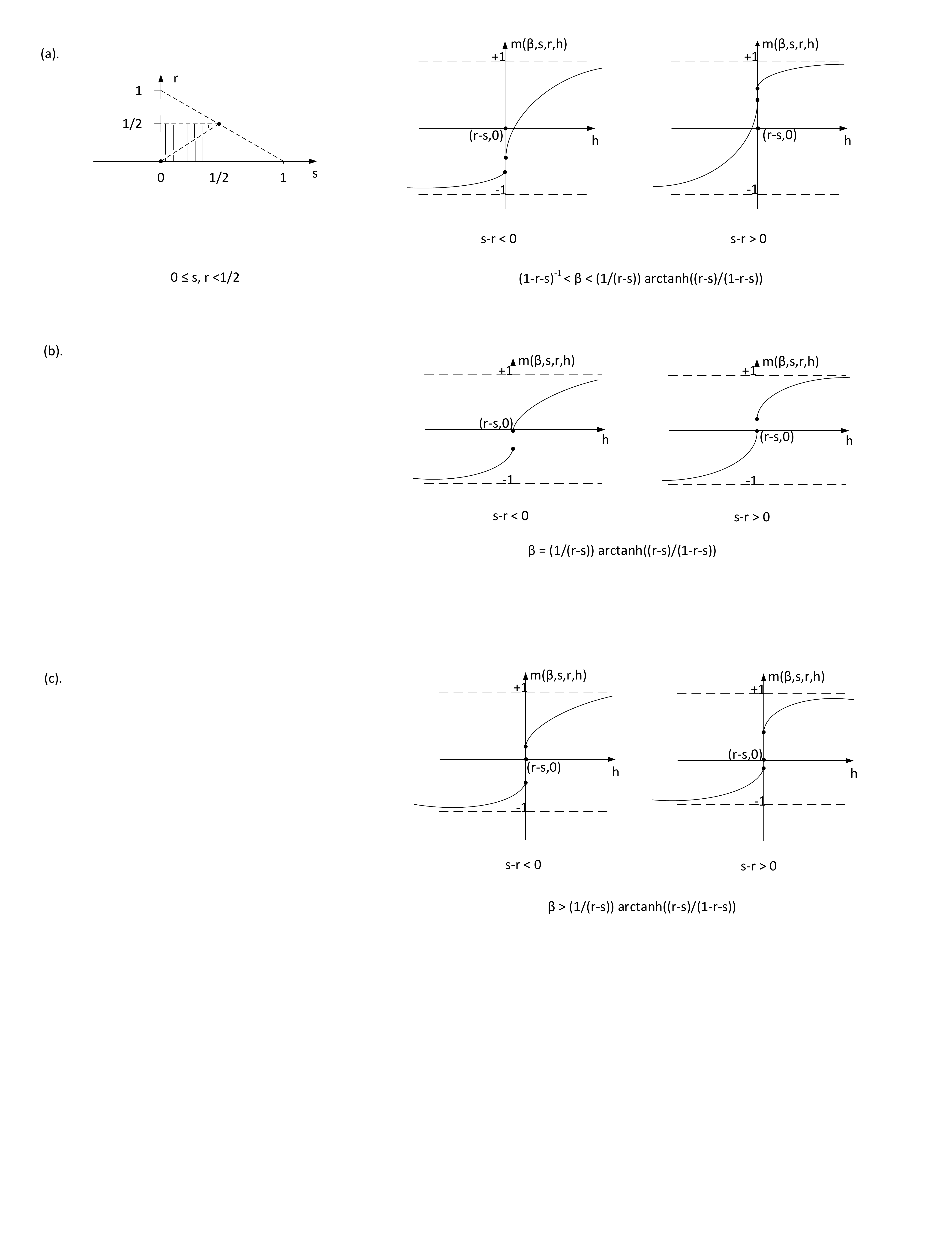}
\end{flushleft}
\vspace{-6.9cm}
\caption{\small \emph{ Discontinuity in the map $h\mapsto m(\beta,s,r,h)$ at $h=r-s$ for $\beta> (1-s-r)^{-1}$
 and for  choices of  $(s,r)$ from the shaded region in the $sr$-plane. The  discontinuity in (a) does not lead 
 to  change in phase/behaviour. That in (b) results in transition from an ordered phase, with positive/ negative 
 magnetization, to disordered phase with zero magnetization.  The situation in (c) leads to transtion from positive 
 magnetization to negative one or vice versa though  $m(\beta,s,r,(r-s)^{-})$ and $m(\beta,s,r,(r-s)^{+})$ 
 are not symmetric about zero.}}
\label{fig4}
\vspace{0cm}
\end{figure}
\section{Discussion}\label{S2}
\begin{enumerate}
\item As we will see in Section \ref{proof-of-mt}, the conditional Curie-Weiss model $\mu_{N,s,r}$ 
is a Curie-Weiss model on the sites in $V_{N,t_N}$ with $N$-dependent coupling strength 
$1-s_N-r_N$ and $N$-dependent external field $h+ s_N-r_N$. Theorem \ref{maintheorem} 
says that, up to  a constant term $s-r$ and a scaling factor  $1-s-r$, the large $N$-limit 
of the magnetization under $\mu_{N,s,r}$ is determined by the magnetization of  a Curie-Weiss 
model with coupling  strength $1-s-r$ and external field $h+s-r$. The scaling factor $1-s-r$ 
is the effective size of the set of sites in $V_{N,t_N}$. Notice that if $s=r$=0,  then we get the original 
Curie-Weiss model with unit coupling strength and external field $h$. 
\item If $s=r\neq \frac12$, then the $s-r$ part of the conditional magnetization vanishes. The 
conditional magnetization is then  equal to  the magnetization of the Curie-Weiss model  with
external field $h$ and coupling strength scaled  by the effective  size $1-2s$ of the
conditional model. Corollary \ref{cor1} says that  in the regime where $s=r\neq \frac12$,  
the singularity  in the magnetization as $h$ goes to zero and  at $\beta$-values 
above $(1-2s)^{-1}$ is similar to that of the original model only that in the conditional 
model the jump in the magnetization is suppressed by a factor of $1-2s$.

\item In the regime where either  $\frac12\leq s<1$ or $\frac12\leq r<1$, with $s+r<1$, it is clear that 
the  conditional model will always be  negatively or positively magnetized depending on the sign of 
$s-r$. Though the conditional model  is always  magnetized along the sign  of $s-r$, yet Corollary \ref{cor2} 
indicates that  the conditional magnetization  discontinuously jumps as $h$ tends to $r-s$ at  
$\beta$-values greater than $(1-s-r)^{-1}$ (see Figure \ref{fig3}). This discontinous jump does not lead 
to phase change.

\item The most interesting region in the $sr$-plane is where  $0\leq s, r<\frac12$ and $s\neq r$.
 In this region, the conditional magnetization has three different forms 
of discontinous jumps depending on the choice  of $\beta$-value(s). Here there is a second 
transition point for $\beta$ namely; 
\begin{equation}\label{betass}
\beta^{**}= \dfrac{1}{r-s}\arctanh\left(\dfrac{r-s}{1-s-r}\right)
\end{equation}
in addition to the earlier transition point 
\begin{equation}\label{betas}
\beta^*=(1-s-r)^{-1}.
\end{equation}
Note that $\beta^*<\beta^{**}$, for  $0\leq s, r<\frac12$ and $s\neq r$. Corollary \ref{cor3} 
says the following:
\begin{enumerate}
  \item For $\beta\in (\beta^*, \beta^{**})$, the discontinuity in the conditional magnetization is 
  analogous to the case in Corollary \ref{cor2}. The conditional model is magnetized along 
  the sign of $s-r$, with jump discontinuity at $h=r-s$. We have a discontinous jump in an order 
  parameter that is not accompanied by change in  phase (see  Figure \ref{fig4}(a)).
  \item In the case $\beta=\beta^{**}$,  there is also discontinuous jump in the magnetization. Here 
  there is a nonzero magnetization as $h$ approaches $r-s$ from the values of $h$ for which $h+s-r$ has the 
  same sign as $s-r$. On the other hand, there is no net magnetization as $h$ goes to $r-s$ from 
  $h$-values for which  $s-r$ and $h+s-r$ have different signs.  Thus the discontinuity here leads 
  to a change from  an ordered phase, with a net magnetization, to a disordered phase, with no 
  magnetization or vice versa (see Figure \ref{fig4}(b)).
  \item Further, the discontinous jump for the case $\beta>\beta^{**}$ is similar to the case with $s=r$ 
  discussed in Corollary \ref{cor1}. This discontinuity is followed by change in phase,  either from negatively 
  magnetized phase to a positively magnetized phase or vice versa. Observe that in  Corollary \ref{cor1} the net 
  magnetizations at the discontinuity point are symmetric about zero. Over here, the net magnetizations 
  at the discontinuity point  are not symmetric about zero, due to the presence of the  term $s-r$ (see Figure \ref{fig4}(c)).

  \end{enumerate}
\end{enumerate}

\section{Proofs}\label{S3}

\subsection{Proof of Theorem \ref{maintheorem}}\label{proof-of-mt}

\begin{proof}
Note from \eqref{configsr}
that for any $\sigma\in\Omega_{N,s,r}$
\begin{equation}\label{mN}
m_{N,s,r}=\frac1N\sum_{i=1}^N\sigma_i=s_N-r_N+\frac1N\sum_{i\in V_{N,t_N}}\sigma_i
\end{equation}
Therefore 
\begin{equation}\label{mN2}
m_{N,s,r}^2=s_N^2+r_N^2-2s_Nr_N+\left(\frac1N\sum_{i\in V_{N,t_N}}\sigma_i\right)^2+
2s_N \frac1N\sum_{i\in V_{N,t_N}}\sigma_i- 2r_N \frac1N\sum_{i\in V_{N,t_N}}\sigma_i.
\end{equation}
It follows from the last equality of \eqref{Hamil}, \eqref{mN} and \eqref{mN2} that for any
 $\sigma\in\Omega_{N,s,r}$, the Hamiltonian $H_N(\sigma)$  can be written as follows:
\begin{equation}\label{Hamilnew}
H_N(\sigma)=- N\left[\frac12 (s_N^2+r_N^2)-s_Nr_N+ h(s_N-r_N)+\frac N2\right]- H_{N,s,r}(\sigma)
\end{equation}

where 
\begin{equation}\label{vntn}
\begin{split}
H_{N,s,r}(\sigma)&=-N\left[\frac{1}{2}\left(\frac1N\sum_{i\in V_{N,t_N}}\sigma_i\right)^2+
\frac{s_N-r_N+h}{N}\sum_{i\in V_{N,t_N}}\sigma_i\right]\cr
&=\frac{|V_{N,t_N}|}{t_N}\left[\frac{1}{2}\left(\frac{t_N}{|V_{N,t_N}|}\sum_{i\in V_{N,t_N}}\sigma_i\right)^2+
\frac{t_N[s_N-r_N+h]}{|V_{N,t_N}|}\sum_{i\in V_{N,t_N}}\sigma_i \right]\cr
&=|V_{N,t_N}|\left[\frac{t_N}{2}\left(\frac{1}{|V_{N,t_N}|}\sum_{i\in V_{N,t_N}}\sigma_i\right)^2+
\frac{s_N-r_N+h}{|V_{N,t_N}|}\sum_{i\in V_{N,t_N}}\sigma_i \right].
\end{split}
\end{equation}
The second equality of \eqref{vntn} follows from our assumption that $t_N=\frac{|V_{N,t_N}|}{N}$. Comparing equality 
three of \eqref{vntn} with the second equality of  \eqref{Hamil} we observe that $H_{N,s,r}$ has an $N$-dependent 
coupling strength of $t_N$ as against $1$ for $H_N$. Further,  $H_{N,s,r}$ has an $N$-dependent 
external field of $h+s_N-r_N$ as against $h$ for $H_N$. 
Therefore substituting this form of $H_N$ in \eqref{Hamilnew} into the expression for $\mu_{N,s,r}$ in \eqref{condCWmodel}
leads to 
\begin{equation}\label{condCWmodelnew}
\begin{split}
\mu_{N,s,r}(\sigma)&=\frac{e^{-\beta H_{N,s,r}(\sigma)}}{\sum_{\eta\in\Omega_{N,s,r}} e^{-\beta H_{N,s,r}(\eta)}}\cr
&=\frac{e^{-\beta H_{N,s,r}(\sigma)}}{Z_{N,s,r}}.
\end{split}
\end{equation}
Thus $\mu_{N,s,r}$ is the Curie-Weiss model on the set of sites in $V_{N,t_N}$ with an $N$-dependent 
external field $h+s_N-r_N$ and $N$-dependent coupling strength $t_N$. In view of this,  we have for 
$\beta>0$, $h+s-r\neq0$ and $\beta\leq (1-s-r)^{-1}$,  $h+s-r=0$ that

\begin{equation}
\begin{split}
m(\beta,s,r,h)&=\lim_{N\rightarrow\infty}\frac1N\int \left(\sum_{i\in V_N}\sigma_i\right)d\mu_{N,s,r}(\sigma)\cr
&= \lim_{N\rightarrow\infty}\left[s_N-r_N+\frac1N\int\left(\sum_{i\in V_{N,t_N}}\sigma_i\right)d\mu_{N,s,r}(\sigma) 
\right]\cr
&=s-r+\lim_{N\rightarrow\infty}\frac{|V_{N,t_N}|}{N} \frac{1}{|V_{N,t_N}|}\int\left( \sum_{i\in V_{N,t_N}}\sigma_i
\right)d\mu_{N,s,r}(\sigma)\cr
&=s-r +(1-s-r)\lim_{N\rightarrow\infty}\frac{1}{|V_{N,t_N}|}\int\left( \sum_{i\in V_{N,t_N}}\sigma_i
\right)d\mu_{N,s,r}(\sigma)\cr
&=s-r +(1-s-r)z(\beta,s,r,h),
\end{split}
\end{equation}
where according to Theorem \ref{mainfact} and the last eqality of \eqref{vntn}, $z(\beta,s,r,h)$ 
is the $z$-value in $[-1,1]$ that minimizes  
\begin{equation}
i_{\beta,s,r,h}(z)=-\frac12\beta(1-s-r)z^2-(h+s-r)z+\frac{1+z}{2}\log(1+z)+\frac{1-z}{2}\log(1-z).
\end{equation}  
Next observe that for $\beta\leq (1-s-r)^{-1} $,  $z(\beta,s,r,h)$ tends to zero as $h$ goes to $r-s$ and 
for $\beta>(1-s-r)^{-1}$, $z(\beta,s,r,h)\mapsto z(\beta,s,r,(r-s)^{+})>0$ as $h\downarrow r-s$ and 
$z(\beta,s,r,h)\mapsto z(\beta,s,r,(r-s)^{-})<0$ as $h\uparrow r-s$.

\end{proof}

\subsection{Proof of Corollary \ref{cor1}}\label{proof-of-cor1}

\begin{proof}
Note that for $s=r$, $m(\beta,s,s,h)=(1-2s)z(\beta,s,s,h)$ and for $\beta>(1-2s)^{-1}$ and $h=0$, 
$i_{\beta,s,s,0}(z)$ has two minimizers $z(\beta,s,s,0^+)>0$ and $z(\beta,s,s,0^-)=-z(\beta,s,s,0^+)$. Further, as 
$h\downarrow0$,  $z(\beta,s,s,h)\rightarrow z(\beta,s,s,0^+)$ and  $z(\beta,s,s,h)\rightarrow z(\beta,s,s,0^-)$ 
as $h\uparrow 0$.

\end{proof}

\subsection{Proof of Corollary \ref{cor2}}\label{proof-of-cor2}

\begin{proof}
As indicated in the proof of Corollary \ref{cor1}, the discontinuity in $m(\beta,s,r,h)$ is 
created by the discontinuity in the  $z(\beta,s,r,h)$ part of $m(\beta,s,r,h)$.  Now suppose 
we are in the regime where $m(\beta,s,r,(r-s)^{\pm})>0$ no matter the choice of $r$, such that 
$s+r<1$. Then the question is what range of $s$-values will permit this behaviour. First of all note that 
the largest possible jump that can occur in this case is when  $z(\beta,s,r,(r-s)^{-})=-1$. If at this value of 
$z(\beta,s,r,(r-s)^{-})$, $m(\beta,s,r,(r-s)^{-})>0$ , then we have that   $s-r-(1-s-r)>0$, which leads to $s>\frac12$. 
Thus if $s>\frac12$, then the minimum value of  $z(\beta,s,r,(r-s)^{-})$ is incapable of making  $m(\beta,s,r,(r-s)^{-})$ 
nonpositive.
Futher, for $s=\frac12$, $$m(\beta,s,r,(r-s)^-)=\left(\frac12-r\right)\left[1+z(\beta,s,r,(r-s)^-)\right]
>0,$$ for $\beta\in (0,\infty)$ as $r<\frac12$  and $z(\beta,s,r,(r-s)^-)>-1$, for $\beta\in(0,\infty)$. 
Thus \\$m(\beta,s,r,(r-s)^{\pm})>0$ for $s\geq \frac12$ and $\beta>0$, though  
$m(\beta,s,r,(r-s)^+)>m(\beta,s,r,(r-s)^-).$ Similarly, $m(\beta,s,r,(r-s)^{\pm})<0$,  for $r\geq\frac12$ 
and $\beta>0$.

\end{proof}

\subsection{Proof of Corollary \ref{cor3}}\label{proof-of-cor3}

\begin{proof}
Observe that the  $|m(\beta,s,r,(r-s)^{\pm})|$ is an increasing function of  $\beta$, with codomain $[-1,1]$.
Suppose  $0\leq r<s<\frac12$. The question now is at what value of $\beta$ is 
$m(\beta,s,r,(r-s)^{-})=0$, thus 
\begin{equation}\label{sreq}
\begin{split}
z(\beta,s,r,(r-s)^-)=\frac{r-s}{1-s-r}.
\end{split}
\end{equation}
Further, $z(\beta,s,r,(r-s)^-)$ is one of the minimizers of the function $z\mapsto i_{\beta,s,r,r-s}(z)$
given by
\[i_{\beta,s,r,r-s}(z)=-\frac{\beta}{2}(1-s-r)z^2+\frac{1+z}{2}\log(1+z)+\frac{1-z}{2}\log(1-z), \quad 
\text{ for }\quad z\in[-1,1].\]
These minimizers satisfy the self-consistency equation 
\begin{equation}\label{zmin}
z=\tanh(\beta [1-s-r]z).
\end{equation}
Therefore, it follows from \eqref{sreq}  and \eqref{zmin}  that the minimizer $z$ 
for which   $m(\beta,s,r,(r-s)^{-})=0$ must  satisfy 
\[ \frac{r-s}{1-s-r}=\tanh\left(\beta(r-s) \right). \]
This  gives rise to the following expression for $\beta$:
\begin{equation}
\beta=\dfrac{1}{r-s}\arctanh\left( \dfrac{r-s}{1-s-r}\right).
\end{equation}
Since $\beta\mapsto |m(\beta,s,r,(r-s^-))|$ is increasing, we have the following:
\begin{enumerate}
\item For $\dfrac{1}{1-s-r}<\beta<\dfrac{1}{r-s}\arctanh\left(\dfrac{r-s}{1-s-r}\right)$,
there is a discontinuity in the map $h\mapsto m(\beta,s,r,h)$, at $h=r-s$, though 
both $m(\beta,s,r,(r-s)^-)$ and  $m(\beta,s,r,(r-s)^+)$ will be positive.
\item At $\beta=\dfrac{1}{r-s}\arctanh\left(\dfrac{r-s}{1-s-r}\right)$, we experience discontiuity 
at $h=r-s$ but this time round $m(\beta,s,r,(r-s)^-)=0$ and  $m(\beta,s,r,(r-s)^+)>0$.
\item  Similar result holds  for $\beta>\dfrac{1}{r-s}\arctanh\left(\dfrac{r-s}{1-s-r}\right)$, with 
$m(\beta,s,r,(r-s)^-)<0$ and  $m(\beta,s,r,(r-s)^+)>0$.
\end{enumerate}
The proof for the case $0\leq s<r<\frac12$ follows  from arguments similar to the case considered 
above.

\end{proof}



\end{document}